\documentclass[12pt]{article}
\font\gordas = msbm10 at 12pt
\def\bbb#1{\hbox {{\gordas #1}}}
\def\a{{\bbb A}}

\def\erre{{\bbb R}}
\def\ce{{\bbb C}}

\def\o{{\bbb O}}
\def\ze{{\bbb Z}}
\def\ache{{\bbb H}}

\usepackage{amssymb}
\begin{document}
\begin{center}{\large\bf Monomorphisms between Cayley-Dickson Algebras}\\[1cm]
Guillermo Moreno\\
Departamento de Matem\'aticas\\
CINVESTAV del IPN\\
M\'exico, D.F.\\
gmoreno@math.cinvestav.mx
\end{center}
\vglue2cm
\noindent
{\bf Abstract}: In this paper we study the algebra monomorphisms from \\
${\a}_m={\erre}^{2^m}$ into ${\a}_n={\erre}^{2^n}$ for $1\leq m\leq n$, where ${\a}_n$ are the Cayley-Dickson algebras. For $n\geq 4$, we show that there are many types of monomorphisms and we describe them in terms of the zero divisors in ${\a}_n$.
\vglue9cm
\noindent
\newpage
\noindent
{\bf Introduction:} The Cayley-Dickson algebras ${\a}_n$ over the real numbers is an algebra structure on ${\erre}^{2^n}={\a}_n$ for $n\geq 0$.

By definition the Cayley-Dickson algebras (C-D algebras) are given by doubling process of Dickson [D].

For $(a,b)$ and $(x,y)$ in ${\a}_n\times {\a}_n$,  define the product in $\a_{n+1}={\a}_n\times{\a}_n$ as follows:
$$(a,b)\cdot(x,y)=(ax-\overline{y}b,ya+b\overline{x}).$$

So if ${\a}_0={\erre}$ and $\overline{x}= x$ for all $x$ in $\erre$ then 
$\a_1={\ce}$ the complex numbers $\a_2={\ache}$ the quaternion numbers and ${\a}_3={\o}$ the octonion numbers.

As is well known ${\a}_n$ is commutative for $n\leq 1$; associative for $n\leq 2$ and alternative for $n\leq 3$, also ${\a}_n$ is normed for $n\leq 3$.

For $n\geq 4\quad  {\a}_n$ is flexible and has zero divisors [Mo$_1$].

Let $Aut({\a}_n)$ be the automorphism group of the algebra ${\a}_n$. 

As is well known 

$Aut({\a}_1)={\ze}/2=\{$ Identity, Conjugation$\}$.

$Aut({\a}_2)=SO(3)$ the rotation group in ${\erre}^3$

$Aut({\a}_3)=G_2$ the exceptional Lie group.

\noindent
(See [H-Y] and [Wh]).

For $n\geq 4$, Eakin-Sathaye showed  that

\hglue3cm $Aut({\a}_n)=Aut({\a}_{n-1})\times \sum_3.$

\noindent
Where $\sum_3$ is the symmetric group of order 6.

(See [E-K] and [Mo$_2$]).

In this paper we will extend the above results in the following sense:

Suppose that $1\leq m\leq n$.

By definition an algebra monomorphism $\varphi:{\a}_m\rightarrow{\a}_n$ is a linear monomorphism such that (i) $\varphi(e_0)=e_0$ and (ii) $\varphi(xy)=\varphi(x)\varphi(y)$ for all $x$ and $y$ in ${\a}_m$ where $e_0=(1,0,\ldots,0)$ is the unit element in ${\a}_m$ and ${\a}_n$ respectively.

We will describe the set
$${\cal M}({\a}_m,{\a}_n)=\{\varphi:{\a}_m\rightarrow{\a}_n|\varphi\quad\hbox{\rm algebra monomorphism}\}.$$
We will see that  this set is more complicated to describe for $n\geq 4$.

For $n\leq 3$ we will recover  known results about the relationship between $Aut({\a}_n)$ and the Stiefel manifolds $V_{2^n-1,2}$.

For $n\geq 3$, recall that $\{e_0,e_1,\ldots,e_{2^n-1}\}$ denotes the canonical basis in ${\a}_n$ and that the doubling process is given by 
$${\a}_{n+1}={\a}_n\oplus{\a}_n\widetilde{e}_0$$ 
where $\tilde{e}_0:=e_{2^n}$ (half of the way basic in ${\a}_{n+1}$).

For $\varphi\in {\cal M}({\a}_m;{\a}_{n+1})$ for $n\geq 3,$

$\varphi$ is of type I 
if $e_{2^n}=\tilde{e}_0\in(Im\varphi)\subset{\a}_{n+1}$ and $\varphi$ is of type II if also $e_{2^n-1}:=\varepsilon\in(Im\varphi)\subset{\a}_{n+1}$.

The main result of this paper is Theorem 2.5: the set of type II monomorphisms
from ${\a}_3$ to ${\a}_{n+1}$ can be described by the set of zero divisors in 
${\a}_{n+1}$ for $n\geq 4$.
\vglue1cm
\noindent
{\bf \S 1. Pure and doubly pure elements in ${\a}_{n+1}.$} 
\vglue.5cm
\noindent

    Throughout this paper we will stablish the following notational conventions:Elements in ${\a}_n$ will be denoted by Latin characters  $a,b,c,\ldots, x,y,z.$
Elements in ${\a}_{n+1}$ will be denoted by Greek characters  $\alpha,\beta,\gamma,\ldots$ For example, 
$$\alpha=(a,b)\in{\a}_n\times{\a}_n.$$

When we need to represent  elements  in ${\a}_n$ as  elements  in ${\a}_{n-1}\times{\a}_{n-1}$ we use subscripts, for instance, $a=(a_1,a_2),\quad b=(b_1,b_2),$   and so on, with  $a_1,a_2,b_1,b_2$ in ${\a}_{n-1}.$
\vglue.2cm
\noindent

Now  $\{e_0,e_1,\ldots,e_{2^n-1}\}$ denotes the canonical basis in ${\a}_n.$ 
Then by the doubling process
$$\{(e_0,0),(e_1,0),\ldots,(e_{2^n-1},0),(0,e_0),\ldots,(0,e_{2^n-1})\}$$
is the canonical basis in ${\a}_{n+1}={\a}_n\times{\a}_n$. By standard abuse of notation, we denote, also $e_0=(e_0,0),e_1=(e_1,0),\ldots,e_{2^n-1}=(e_{2^n-1},0),\\
 e_{2^n}=(0,e_0),\ldots,e_{2^{n+1}-1}=(0,e_{2^n-1})$ in ${\a}_{n+1}.$ 

\vskip.5cm
\noindent

 For $\alpha=(a,b)\in{\a}_n\times{\a}_n={\a}_{n+1}$ we denote $\widetilde{\alpha}=(-b,a)$ (the complexification of $\alpha$) so $\widetilde{e}_0=(0,e_0)$ and $\alpha\widetilde{e}_0=(a,b)(0,e_0)=(-b,a)=\widetilde{\alpha}$. Notice that $\widetilde{\widetilde{\alpha}}=-\alpha.$

\vskip.5cm
\noindent

\mbox{The {\bf trace}} on ${\a}_{n+1}$ is the linear map $t_{n+1}:{\a}_{n+1}\rightarrow{\erre}$ given by $t_{n+1}(\alpha)=\alpha+\overline{\alpha}=2$(real part of $\alpha$) 
so $t_{n+1}(\alpha)=t_n(a)$ when $\alpha=(a,b)\in{\a}_n\times{\a}_n.$
\vglue.5cm
\noindent
{\bf Definition:} $\alpha=(a,b)$ in ${\a}_{n+1}$ \underline{is pure} if
$$t_{n+1}(\alpha)=t_n(a)=0.$$

$\alpha=(a,b)$ in ${\a}_{n+1}$ \underline{is doubly pure} if it is pure and also $t_n(b)=0$;  i.e., $\widetilde{\alpha}$ is pure  in ${\a}_{n+1}.$

Also $2\langle a,b\rangle =t_n (a\overline{b})$ for $\langle  ,  \rangle $ the inner product in ${\erre}^{2^n}$ (see [Mo$_1$]).

Notice that for $a$ and $b$ pure elements $a\perp b$ if and only if $ab=-ba.$
\vglue.5cm
\noindent
{\bf Notation:} $Im({\a}_n)=\{e_o\}^\perp\subset{\a}_n$ is the vector subspace consisting of pure elements in ${\a}_n$;  i.e., 
$Im({\a}_n)={\rm Ker}(t_n)={\erre}^{2^n-1}.$

$\widetilde{{\a}}_{n+1}=Im({\a}_n)\times Im({\a}_n)=\{e_0,\widetilde{e}_0\}^\perp={\erre}^{2^{n+1}-2}$ is the vector subspace consisting of doubly pure elements in ${\a}_{n+1}.$
\vglue.3cm
\noindent
{\bf Lemma 1.1.} For $a$ and $b$ in $\widetilde{{\a}_n}$ we have that
\begin{itemize}
\item[1)] $a\widetilde{e}_0=\widetilde{a}$ and $\widetilde{e}_0a=-\widetilde{a}.$
\item[2)] $a\widetilde{a}=-||a||^2\widetilde{e}_0$ and $\widetilde{a} a=||a||^2\widetilde{e}_0$ so $a\perp \widetilde{a}.$
\item[3)] $\widetilde{a}b=-\widetilde{ab}$ with $a$ a pure element.
\item[4)] $a\perp b$ if and only if $\widetilde{a}b+\widetilde{b}a=0.$
\item[5)] $\widetilde{a}\perp b$ if and only if $ab=\widetilde{b}\widetilde{a}.$
\item[6)] $\widetilde{a}b=a\widetilde{b}$ if and only if $a\perp b$ and $\widetilde{a}\perp b.$
\end{itemize}
{\bf Proof:} Notice that $a$ is pure if $\overline{a}=-a$ and if $a=(a_1, a_2)$ is doubly pure, then $\overline{a}_1=-a_1$ and $\overline{a}_2=-a_2.$
\begin{itemize}
\item[1)] $\widetilde{e}_0a=(0, e_0)(a_1, a_2)=(-\overline{a}_2, \overline{a}_1)=(a_2, -a_1)=-(-a_2,a_1)=-\widetilde{a}.$
\item[2)] $a\widetilde{a}=(a_1, a_2) (-a_2, a_1)=(-a_1a_2+a_1a_2, a_1^2+a^2_2)=(0,-||a||^2e_0)=-||a||^2\widetilde{e}_0.$

Similarly $\widetilde{a}a=(-a_2, a_1)(a_1, a_2)=(-a_2a_1+a_2a_1, -a^2_2-a^2_1)=||a||^2\widetilde{e}_0.$ 

Now, since $-2\langle \widetilde{a}, a\rangle =a\widetilde{a}+\widetilde{a}a=0$ we have $a\perp \widetilde{a}.$
\item[3)] $\widetilde{a}b=(-a_2, a_1)(b_1, b_2)=(-a_2b_1+b_2a_1, -b_2a_2 -a_1b_1)$.

 So $\widetilde{\widetilde{a}b}= (a_1b_1+b_2a_2, b_2a_1-a_2b_1)=(a_1, a_2)(b_1, b_2)=ab$ then 
$-\widetilde{a}b=\widetilde{ab}.$

Notice that in this proof we only use that $\overline{a}_1=-a_1$; i.e., $a$ is pure and $b$ doubly pure.
\item[4)] $a\perp b\Leftrightarrow ab+ba=0\Leftrightarrow ab=-ba\Leftrightarrow\widetilde{ab}=-\widetilde{ba}.$

$\Leftrightarrow -\widetilde{a}b=\widetilde{b}a\Leftrightarrow \widetilde{a}b+\widetilde{b}a=0$ by (3).
\item[5)] $\widetilde{a}\perp b \Leftrightarrow \widetilde{\widetilde{a}}b+\widetilde{b}\widetilde{a}=0$ (by (4)) $\Leftrightarrow -ab+\widetilde{b}\widetilde{a}=0.$
\item[6)] If $\widetilde{a}\perp b$ and $a\perp b$, then by (3) and (4)
$
\widetilde{a}b=-\widetilde{ab}=\widetilde{ba}=-\widetilde{b}a=a\widetilde{b}.
$

Conversely, put $a=(a_1, a_2)$ and $b=(b_1, b_2)$ in ${\a}_{n-1}\times
{\a}_{n-1}$ and define 

$c:=(a_1b_1+b_2a_2)$ and $d:=(b_2a_1-a_2b_1)$ in ${\a}_{n-1}.$

Then $a\widetilde{b}=(a_1, a_2)(-b_2, b_1)=(-a_1b_2+b_1a_2, b_1a_1+a_2b_2)$ so $a\widetilde{b}=(-\overline{d}, \overline{c}).$
\end{itemize}
Now $ab=(a_1, a_2)(b_1, b_2)=(a_1b_1+b_2a_2, b_2a_1-a_2b_1)=(c,d)$, so $\widetilde{ab}=(-d,c)$ and then $\widetilde{a}b=(d, -c)$. 

Thus, if $a\widetilde{b}=\widetilde{a}b$ then $\overline{c}=-c$ and $d=-\overline{d}.$
Then 
\begin{eqnarray*} 
t_n(ab)&=&t_{n-1}(c)=c+\overline{c}=0\;\;\hbox{\rm and}\;\;a\perp b\\
t_n(\widetilde{a}b)&=&t_{n-1}(d)=d+\overline{d}=0\;\;\hbox{\rm and}\;\;\widetilde{a}\perp b.
\end{eqnarray*}
\hfill Q.E.D

\vskip .3cm

\noindent
{\bf Corollary 1.2} For each $a\neq 0$ in $\widetilde{{\a}}_n$. The four dimensional vector subspace generated by $\{e_0, \widetilde{a}, a, \widetilde{e}_0\}$ is a copy of ${\a}_2 ={\ache}.$(We denote it by ${\ache}_a).$
\vglue .25cm
\noindent
{\bf Proof:} We suppose that $||a||=1$, otherwise we take ${a\over ||a||}$. Construct the following multiplication table.

\vskip .5cm
\begin{center}
\begin{tabular}{c|cccc}
 & $e_0$ & $\widetilde{a}$ & $a$ & $\widetilde{e}_0$\\ \hline
$e_0$ & $e_0$ & $\widetilde{a}$ & $a$ & $\widetilde{e}_0$ \\
$\widetilde{a}$ & $\widetilde{a}$ & $-e_0$ & $+\widetilde{e}_0$ & $-a$\\
$a$ & $a$ & $-\widetilde{e}_0$ & $-e_0$ & $\widetilde{a}$\\
$\widetilde{e}_0$ & $a$ & $\widetilde{e}_0$ & $-\widetilde{a}$ & $-e_0$
\end{tabular}
\end{center}

By lemma 2.1. $a\widetilde{e}_0=\widetilde{a}$; $\widetilde{e}_0a=-\widetilde{a}$; $\widetilde{a}\widetilde{e}_0=\widetilde{\widetilde{a}}=-a;$ $\widetilde{e}_0 \widetilde{a}=-\widetilde{\widetilde{a}}=a$; $a\widetilde{a}=-\widetilde{e}_0$ and $\widetilde{a}a=\widetilde{e}_0.$

But this is the multiplication table of ${\a}_2 ={\ache}$ identifying $e_0\leftrightarrow 1$, $\widetilde{a}\leftrightarrow \hat{i},$ $a\leftrightarrow \hat{j}$ and $\widetilde{e}_0\leftrightarrow \hat{k}.$

\hfill Q.E.D.

\vskip.5cm
\noindent
{\bf  $\S$ 2. Monomorphism from ${\a}_m$ to ${\a}_n$.}
\vglue.3cm
Throughout this chapter $1\leq m\leq n$.
\vglue.3cm
\noindent
{\bf Definition.} An algebra monomorphism from 
${\a}_m$ to ${\a}_n$ is a linear monomorphism 
$\varphi:{\a}_m\rightarrow{\a}_n$ such that 

i) $\varphi(e_0)=e_0$ (the first $e_0$ is in ${\a}_m$ and the second 
$e_0$  in ${\a}_n$)

ii) $\varphi(xy)=\varphi(x)\varphi(y)$ for all $x$ and $y$ in ${\a}_m$.

By definition we have that $\varphi(re_0)=r\varphi(e_0)$ for all $r$ in 
${\erre}$ so \\
$\varphi (Im({\a}_m))\subset\varphi (Im({\a}_n))$ and $\varphi(\overline{x})=\overline{\varphi(x)}$.

Therefore $||\varphi(x)||^2=\varphi(x)\overline{\varphi(x)}=\varphi(x)\varphi(\overline{x})=\varphi(x\overline{x})=\varphi(||x||^2)=||x||^2$ for all 
$x\in {\a}_m$ and $||\varphi(x)||=||x||$ and $\varphi$ is an orthogonal linear transformation from ${\erre}^{2^m-1}$ to ${\erre}^{2^n-1}$.

The \underline{trivial} monomorphism is the one given by $\varphi(x)=
(x,0,0,\ldots,0)$ for $x\in{\a}_m$ and $0$ in 
${\a}_m \quad (2^{n-m-1}$-times$)$.

$ {\cal M}({\a}_m,{\a}_n)$ denotes the set of algebra monomorphisms from ${\a}_m$ to ${\a}_n$.

For 
$m=n,\,\,{\cal M}({\a}_m;{\a}_n)={\rm Aut}({\a}_n)$ 
the group of algebra automorphisms of ${\a}_n$
\vglue.3cm
\noindent
{\bf Proposition 2.1.} ${\cal M}({\a}_1;{\a}_n)=S(Im(
{\a}_n))=S^{2^n-2}$.
\vskip.3cm
\noindent
{\bf Proof:} ${\a}_1={\ce}=Span\{e_0,e_1\}$.

If $x\in{\a}_1$ then $x=re_0+se_1$ and for $w\in Im({\a}_n)$ with $||w||=1$ we have that $\varphi_w(x)=re_0+sw$ define an algebra monomorphism from ${\a}_1$ to ${\a}_n$. This can be seen by direct calculations, recalling that, Center $({\a}_n)={\erre}$ for all $n$ and that every associator with one real entrie vanish.
\begin{eqnarray*}
\varphi_w(x)\varphi_w(y)&=&(re_0+sw)(pe_0+qw)=(rp+sqw^2)e_0+(rq+sp)w\\
&=&(rp-sq)e_0+(rq+sp)w\\
&=&\varphi_w(x)\varphi_w(y)
\end{eqnarray*}
when $y=pe_0+qe_1$ and $p$ and $q$ in ${\erre}$. Clearly $\varphi_w(e_0)=e_0$. 

Conversely, for $\varphi\in{\cal M}({\a}_1;{\a}_n)$, set $w=\varphi(e_1)$ in ${\a}_n$ so $||w||=1$ and $\varphi_w=\varphi$.

\hfill Q.E.D.

\vskip.5cm
\noindent
{\bf Remark:} In particular, we have that
$$Aut({\a}_1)=S^0={\ze}/2=\{\hbox{\rm Identity, conjugation}\}=\{\varphi_{e_1},\varphi_{-e_1}\}.$$
To calculate ${\cal M}({\a}_2;{\a}_n)$ for $n\geq 2$ we need to recall (see [Mo$_2$]).
\vglue.5cm
\noindent
{\bf Definition:} For $a$ and $b$ in ${\a}_n$. We said that 
\underline{$a$ alternate with $b$}, we denote it by $a\rightsquigarrow b$, if $(a,a,b)=0$.

We said that \underline{$a$ alternate strongly with $b$}, we denote it by $a\leftrightsquigarrow b$, if $(a,a,b)=0$ and $(a,b,b)=0$.

Clearly $a$ alternate strongly with $e_0$ for all $a$ in ${\a}_n$ and if $a$ and $b$ are linearly dependent then $a\leftrightsquigarrow b$ (by flexibility).

Also,by definition,$a$ is an alternative element if and only if $a\rightsquigarrow x$ for all $x$ in ${\a}_n$.

By Lemma 1.1 (1) and (2) we have that for any doubly pure element $a$ in 
${\a}_n\,\, (a,a,\tilde{e}_0)=0$ and (by the above remarks) 
$\tilde{e}_0$ alternate strongly with any $a$ in ${\a}_n$.

For $a$ and $b$ pure elements in ${\a}_n,$ we define the vector subspace of ${\a}_n$
$$V(a;b)={\rm Span}\{e_0,a,b,ab\}.$$
Also we identify the Stiefel manifold $V_{2^n-1,2}$ as 
$$\{(a,b)\in Im({\a}_n)\times Im({\a}_n)|a\perp b, ||a||=||b||=1\}$$
{\bf Lemma 2.2.} If $(a,b)\in V_{2^n-1,2}$ and $a\leftrightsquigarrow b$ then $V(a;b)={\a}_2={\ache}$ the quaternions.

\vglue.5cm
\noindent
{\bf Proof:} Suppose that $(a,b)\in V_{2^n-1,2}$ and that $(a,a,b)=0$ then we have
\begin{eqnarray*}
\langle ab,a\rangle&=&\langle b,\overline{a}a\rangle=\langle b,||a||^2e_0\rangle=||a||^2\langle b,e_0\rangle=0\\
\langle ab,a\rangle&=&\langle a,b\overline{b}\rangle=\langle a,||b||^2e_0\rangle=||b||^2\langle a,e_0\rangle=0\\
||ab||^2=\langle ab,ab\rangle&=&\langle\overline{a}(ab),
b\rangle=\langle-a(ab),b\rangle=\langle-a^2b,b\rangle\\
&=&-a^2\langle b,b\rangle=||a||^2||b||^2=1
\end{eqnarray*}
so $\{e_0,a,b,ab\}$ is an orthonormal set of vectors in ${\a}_n$.

Finally using also that $(a,b,b)=0$ and $ab=-ba$ we may check by direct calculations that the multiplication table of $\{e_0,a,b,ab\}$ coincides with the one of the quaternions and by the identification $e_0\mapsto e_0, a\mapsto e_1, b\mapsto 
e_2$ and $ab\mapsto e_3$ we have an algebra isomorphism between 
${\a}_2={\ache} $ and $V(a;b)$.

\hfill Q.E. D.
\vglue.5cm
\noindent
{\bf Proposition 2.3.} ${\cal M}({\a}_2;{\a}_n)=\{(a,b)\in V_{2^n-1,2}|a\leftrightsquigarrow b\}$ for $n\geq 2$.

In particular
$$Aut({\a}_2)={\cal M}({\a}_2;{\a}_2)=V_{3,2}=SO(3)$$
and 
$${\cal M}({\a}_2,{\a}_3)=V_{7,2}.$$
{\bf Proof.} The inclusion ``$\supset$'' follows from Lemma 2.2. Conversely suppose that $\varphi \in{\cal M}({\a}_2,{\a}_n)$ then $\varphi(e_0)=e_0,  (\varphi(e_1),\varphi(e_2))\in V_{2^n-1,2}$ and 
$$V(\varphi(e_1),\varphi(e_2))={\rm Im}\varphi={\ache}\subset{\a}_n.$$

Since ${\a}_2$ is an associative algebra and ${\a}_3$ is an alternative algebra we have that $a\leftrightsquigarrow b$ for any two elements in ${\a}_n$ for $n=2$ or $n=3$.

\hfill Q.E.D.
\vglue.3cm
\noindent
{\bf Remark.} Recall that $\tilde{{\a}}_n=\{e_0,\tilde{e}_0\}^\perp={\erre}^{2^n-2}$ denotes the vector subspace of doubly pure elements. Since $a\leftrightsquigarrow
\tilde{e}_0$ for any element in $\tilde{{\a}}_n,$ we have that, if $a\in S(\tilde{{\a}}_n)$ i.e., $||a||=1$ then $(a,\tilde{e}_0)\in V_{2^n-1,2}$ and the assignment $a\mapsto (a,\tilde{e}_0)$ define an inclusion from 
$S(\tilde{{\a}}_n)=S^{2^n-3}\hookrightarrow{\cal M}({\a}_2;{\a}_n)\subset V_{2^n-1,2}$ which resembles ``the bottom cell'' inclusion in $V_{2^n-1,2}$.
\vskip.2cm
To deal with the cases $3=m\leq n$ we have to use the notion of a \underline{special triple} (see [Wh] and [Mo$_1$]).
\vskip.2cm
\noindent
{\bf Definition:} A set $\{a,b,c\}$ in 
$Im({\a}_n)$ is \underline{a special triple} if 

(i) $\{a,b,c\}$ is an orthonormal set

(ii) $a\leftrightsquigarrow b, a\leftrightsquigarrow c$ and $b\leftrightsquigarrow c$ i.e. its elements alternate strongly, pairwise.

(iii) $c\in V(a;b)^\perp\subset {\a}_n$.

Now is easy to see that  if $\{a,b,c\}$ is a special triple then
$V(a;b);V(a,c);V(b,c)$ are isomorphic to ${\a}_2$.

For a special triple $\{a,b,c\}$ consider the following vector subspace of ${\a}_n$
$${\o}(a;b;c):=Span \{e_0,a,b,ab,c(ab),cb,ac,c\}.$$
\vskip.2cm
\noindent
{\bf Proposition 2.4:} For a special triple $\{a,b,c\}$ in ${\a}_n$ and $n\geq 3;
\,\,\, {\o}(a,b,c)$ is an eight-dimensional vector subspace isomorphic, as algebra, to ${\a}_3={\o}$ the octonions and ${\cal M}({\a}_3,{\a}_n)=\{(a,b,c)\in
({\a}_n)^3|\{a,b,c\}$ special triple$\}$.
\vskip.2cm
\noindent
{\bf Proof:}  We know that all elements in $\{e_0,a,b,ab,c(ab),cb,ac,c\}$ are of norm one and also
$$\langle c(ab),a\rangle=-\langle ab,ac\rangle=\langle a(ab),c\rangle=\langle a^2b,c\rangle=a^2\langle b,c\rangle=0.$$
Similarly $(c(ab))\perp b$ and $(c(ab))\perp c$. Thus $\{e_0,a,b,ab,c(ab),cb,ac,c\}$ is an orthonormal set of vectors and ${\o}(a;b;c)$ is eight-dimensional. To see that ${\o}(a;b;c)\cong{\a}_3$ we have to construct the corresponding multiplication table,which is, a routine calculation.(See [Mo$_2$])

Conversely, if $\varphi\in{\cal M}({\a}_3,{\a}_n)$ then $a=\varphi(e_1),b=\varphi(e_2)$ and $c=\varphi(e_7)$ form an special triple, when $\{e_0,e_1,e_2,e_3,e_4,e_5,e_6,e_7\}$ is the canonical basis in ${\a}_3,$ and we recall that $e_1e_2=e_3,e_7e_3=e_4, e_7e_2=e_5$ and $e_1e_7=e_6$ in ${\a}_3$.

\hfill Q.E.D.

\vskip.5cm
\noindent
{\bf Remark.} For $n=3, {\a}_3$ is an alternative algebra so a special triple in  
${\a}_3$ is every triple such that 

(i) $\{a,b,c\}$ is orthonormal  

(ii)
$c\perp (ab)$.

So Proposition 4.4 gives the construction of $G_2=Aut({\a}_3)$ as in [Wh] and the assignment
\begin{eqnarray*}
G_2=Aut({\a}_3)&\stackrel{\pi}{\rightarrow}&{\cal M}({\a}_2,{\a}_3)=V_{7,2}\\
(a,b,c)&\mapsto&(a,b)
\end{eqnarray*}
is the known fibration $G_2\stackrel{\pi}{\rightarrow}V_{7,2}$ with fiber 
$S^3$.
\vglue.5cm
\noindent
{\bf Remark.} Suppose that $n\geq 4$ and that $\{a,b,c\}$ is a special triple in ${\a}_n$ so ${\o}(a;b;c)$ is the image of some algebra monomorphism from ${\a}_3$ to ${\a}_n$ and any orthonormal triple $\{x,y,z\}$ of pure elements in 
${\o}(a;b;c)$ with $z\perp(xy)$ is also a special triple in ${\a}_n$ and 
$${\o}(x;y;z)={\o}(a;b;c).$$
{\bf Definition:} For $1\leq m\leq n.$
 $ \varphi\in{\cal M}({\a}_m,{\a}_n)$ 
\underline{is a type I monomorphism} if $\widetilde{e}_0\in $ (Image of $\varphi$) $\subset{\a A}_n$.

Since $\tilde{e}_0=e_{2^{n-1}}$ in ${\a A}_n$ then the trivial monomorphism is \underline{not} a type I monomorphism unless $n=m$,  because by definition  its image is generated by $\{e_0,e_1,\ldots,e_{2^{m-1}-1}\}.$

Denote by ${\cal M}_1({\a A}_m;{\a A}_n)=\{\varphi\in{\cal M}({\a A}_m;{\a A}_n)|\tilde{e}_0\in({\rm Im}\varphi)\}$ the subset of all type I monomorphisms, clearly ${\cal M}_1({\a A}_n,{\a A}_n)={\cal M}({\a A}_n;{\a A}_n)={\rm Aut}({\a A}_n)$. 

Using proposition 2.1 we may verify that for $n\geq 2$. 
$${\cal M}_1({\a A}_1,{\a A}_n)=\{re_0+s\tilde{e}_0|r^2+s^2=1\}=S^1.$$ 
Also,by Lemma 2.2, for a non-zero $a\in\tilde{{\a A}}_n$ we have that $V(a;\tilde{e}_0)={\a H}_a$ and $${\cal M}_1({\a A}_2;{\a A}_n)=S(\tilde{{\a A}}_n)=S^{2^n-3}$$.
In particular
$$S^5={\cal M}_1({\a A}_2;{\a A}_3)\subset{\cal M}({\a A}_2;{\a A}_3)=V_{7,2}$$
is ``the bottom cell'' of $V_{7,2}$.

Also we can check,using the fact that $a\rightsquigarrow\tilde{e}_0$ for all $a\in\tilde{{\a A}}_n,$ that $${\cal M}_1({\a A}_3;{\a A}_n)={\cal M}({\a A}_2;{\a A}_n)\cap V_{2^n-2,2}$$
where $V_{2^n-2,2}=\{(a,b)\in V_{2^n-1,2}|a$ and $b$ are in $\tilde{{\a A}}_n\}$.
\vglue.5cm
\noindent
{\bf Definition:} A  type I
 monomorphism $\varphi\in{\cal M}_1({\a A}_m;{\a A}_{n+1})$ 
\underline{is of type II} if 
$$e_{2^{n-1}}:=\varepsilon\in {\rm (Im}\quad \varphi)\subset {\a A}_{n+1}.$$

By definition if  $\varphi\in{\cal M}_1({\a A}_m;{\a A}_{n+1})$ then $\tilde{e}_0\in (Im\,\,\varphi)\subset{\a A}_{n+1}$ and if we also assume that $\varepsilon\in (Im \varphi)\subset{\a A}_{n+1}$ then $\varepsilon\tilde{e}_0=\tilde{\varepsilon}$ and 
${\a H}_\varepsilon:=Span\{e_0,\tilde{\varepsilon},\varepsilon,\tilde{e}_0\}$ lies in $(Im\,\varphi)\subset{\a A}_{n+1};$ therefore 
$$\varphi\in{\cal M}({\a A}_m,{\a A}_{n+1})\,\,\hbox{\rm 
is of type II if and only if }\,\,
{\a H}_\varepsilon\subset(Im\,\,\varphi).$$

Denote ${\cal M}_2({\a A}_m,{\a A}_{n+1})=\{\varphi\in{\cal M}({\a A}_m;{\a A}_{n+1})|\varphi$ is type II$\}$.
\vglue.5cm
\noindent
{\bf Theorem 2.5} For $n\geq 3$.
$${\cal M}_2({\a A}_3;{\a A}_{n+1})={\a C} P^{2^{n-1}-1}\cup\overline{X}_n$$
where ${\a C} P^m$ is the complex projective space in ${\a C}^m$ and 
$$\overline{X}_n=\{(x,y)\in{\a A}_n\times{\a A}_n|xy =0,x\neq 0\,\,{\rm  and}\,\,
y\neq 0\}.$$ 
In particular for $n=3\quad\overline{X}_3=\Phi$ (empty set) and 
${\cal M}_2({\a A}_3;{\a A}_4)={\a C} P^3$.
\vglue.5cm
\noindent
{\bf Proof:} Suppose that $\varphi:{\a A}_3\rightarrow{\a A}_{n+1}$ is an algebra monomorphism for $n\geq 3$ with ${\a H}_\varepsilon\subset (Im\,\varphi)$. So $Im\,\varphi$ is isomorphic to ${\a O}={\a A}_3,$ as algebras, then there is a non-zero $\alpha\in{\a H}^\perp_\varepsilon\subset\tilde{A}_{n+1}$ such that
$$Im\,\varphi={\a O}_\alpha:=Span \{e_0,\tilde{\varepsilon},\varepsilon,\tilde{e}_0,\tilde{\alpha},\alpha\varepsilon,\tilde{\varepsilon}\alpha,\alpha\}\subset{\a A}_{n+1}$$

Suppose that $\alpha=(a,b)\in\tilde{{\a A}}_n\times\tilde{
{\a A}}_n={\a H}^\perp_\varepsilon\subset{\a A}_{n+1}$ 
and that $a\neq 0$ (similarly we may assume $b\neq 0$). 
 
Now $Im\,\varphi={\a O}_\alpha={\a A}_3$ if and only if 
$$(\alpha,\alpha,\varepsilon):=\alpha^2\varepsilon-\alpha(\alpha\varepsilon)=0
\quad{\rm 
 i.e.}\quad
-||\alpha||^2\varepsilon=\alpha(\alpha\varepsilon).$$

Using Lemma 1.1 we have that
\begin{eqnarray*}
\alpha(\alpha\varepsilon)&=&(a,b)[(a,b)(\tilde{e}_0,0)]=(a,b)[(\tilde{a},-\tilde{b})]=\\
&=&(a\tilde{a}-\tilde{b}b,-\tilde{b}a-b\tilde{a})=(-||a||^2\tilde{e}_0-||b||^2\tilde{e}_0,0)
-(0,(b\tilde{e}_0)a-b(\tilde{e}_0a))\\
&=&-||\alpha||^2\varepsilon+(0,(a,\tilde{e}_0,b)).
\end{eqnarray*}
Therefore $(a,\tilde{e}_0,b)=0$ in ${\a A}_n$ if and only if $(\alpha,\alpha,\varepsilon)=0$ in ${\a A}_{n+1}$.
\vglue.5cm
Since ${\a A}_n={\a H}_a\oplus{\a H}^\perp_a$ we have that $b=c+d$ where $c\in
{\a H}_a$ and $d\in{\a H}^\perp_a$ with $c$ doubly pure i.e., $c\in Span\{a,\tilde{a}\}\subset{\a H}_a$. 

Since ${\a H}_a$ is associative we have that
\begin{eqnarray*}
(a,\tilde{e}_0,b)&=&(a,\tilde{e}_0,c+d)=(a,\tilde{e}_0,c)+(a,\tilde{e}_0,d)=0+(a,\tilde{e}_0,d)\\
&=&(a,\tilde{e}_0,d)=\tilde{a}d+a\tilde{d}=2\tilde{a}d
\end{eqnarray*}
by Lemma 1.1 (1), (2) and (6).
\vglue.5cm
Therefore $(\alpha,\alpha,\varepsilon)=0$ in ${\a A}_{n+1}$ 
if and only if $ad=\tilde{a}d=0$ in ${\a A}_n$.
\vglue.5cm
So, we have two cases: namely $d=0$ and $d\neq 0$ in ${\a A}_n$.

Suppose $d=0$. Thus $b=c\in{\a H}_a$ and $b\in Span\{a,\tilde{a}\}$ so $(a,b)$ determine a complex line in $\tilde{{\a A}}_n={\erre}^{2^n-2}={\a C}^{2^{n-1}-1}$ and $\alpha\in {\a C} P^{2^{n-1}-1}$. 

Suppose $d\neq 0$ (so $b\neq 0$ and $a\neq 0$).

Thus $ad=0$ and $(a,d)\in\overline{X}_n$.

\hfill Q.E.D.

\newpage
\noindent
{\bf References}

\begin{itemize}
\item[{\rm [D]}]L.E. Dickson. On quaternions and their generalization and the history of the 8 square theorem. Annals of Math; 20 155-171, 1919.
\item[{\rm [E-K]}]Eakin-Sathaye. On automorphisms and derivations of Cayley-Dickson algebras. Journal of Pure and Applied Algebra 129, 263-280, 1990.
\item[{\rm [H-Y]}]S.H. Khalil and P. Yiu. The Cayley-Dickson algebras: A theorem of Hurwitz and quaternions. Bolet\'{\i}n de la Sociedad de L\'odz Vol. XLVIII 117-169. 1997.
\item[{\rm [Mo$_1$]}]G. Moreno. The zero divisors of the Cayley-Dickson algebras over the real numbers. Bol. Soc. Mat. Mex. 4 13-27. 1998
\item[{\rm [Mo$_2$]}]G. Moreno. Alternative elements in the Cayley-Dickson algebras.\\
 hopf.math.purdue.edu/pub/Moreno.2001 and math.RA/0404395.
\item[{\rm [Wh]}] G. Whitehead. Elements of Homotopy theory.\\
 Graduate text in Math. 61 Springer-Verlag.
\end{itemize}
\end{document}